# (10, k)- Reverse Multiple

Madline Al- Tahan

Abstract

We consider the integers having the property of reversing its digits when multiplied by a specific integer $k$.
First, we proved that $k$ should be either 1, 4 or 9. Second, we classify these integers as *(10, 1)- reverse* multiples, *(10, 4)- reverse* multiples and *(10, 9)-* reverse multiples. Then we conclude their general form.

## 1. Introduction

Integers have been of important interest for many researchers in the fields of Computational Mathematics and Number Theory. In which many new definitions of integers with special properties have been introduced such as, Palindrome integers and *(g, k)- reverse multiples* (see [2] and [3]).

A *palindrome integer* is an integer that reads the same forward and backward. The first 20 palindrome integers are the first 10 digits and the numbers 11, 22, 33, 44, 55, 66, 77, 88, 99, 101.

For integers (g, k) with $g \geq 3, 1 \leq k < g$, a number $M$ is called a *(g, k)- reverse multiple* if the reversal of the base-g expansion of $M$ is the base-g expansion of the product $kM$.

A *(10, k)- reverse multiple* is an integer satisfying the property of reversing its digits when multiplied by a positive integer k. The decimal numbers 1089 and 2178 are the two smallest *(10, k)*-reverse multiples, their reversals being respectively $9801 = 9 \times 1089$ and $8712 = 4 \times 2178$. There are no other 4-digit examples in base 10 (see [1] and [4]).

Many researchers work on *(g, k)- reverse multiple*. In particular, Young [6], in 1992, introduced certain trees in order to study the problem of finding all *(g, k)*-reverse multiples for a fixed base g. In 2012, R. Webster and G. Williams [6] found all *(10, 9)*-reverse multiples using the reverse divisors. Then, in 2013, Sloane [4] extends Young's work by replacing the trees with certain finite directed graphs known as "Young graphs". In this paper, we find all *(10, k)*-reverse multiples using an original proof that differs from all previous proofs discussed before. Our result which is stated in our main theorem, Theorem 5, coincides with results in [6].

In section 2, we find sufficient condition for the form of *(10, k)- reverse multiple* integers and values of k which coincides with results in [6]. In section 3, we find necessary condition for an integer to be *(10, k)- reverse multiple*. And conclude with our main theorem, Theorem 5.

## 2. Sufficient conditions

We start this section with a Theorem and some Lemmas in which their proof is straightforward.

**Theorem 1.** An integer is *(10, 1)*- reverse multiple if and only if it is a Palindrome.

**Lemma 1.** If $M$ is an integer of the form

$$21\underbrace{9...9}_{l_1}78\underbrace{0...0}_{m_1}\cdots\underbrace{0...0}_{m_{i-1}}21\underbrace{9...9}_{l_i}78\underbrace{0...0}_{m_i}21\underbrace{9...9}_{l_i}78\underbrace{0...0}_{m_{i-1}}\cdots\underbrace{0...0}_{m_1}21\underbrace{9...9}_{l_1}78$$

or

$$21\underbrace{9...9}_{l_1}78\underbrace{0...0}_{m_1}\cdots 21\underbrace{9...9}_{l_{i-1}}78\underbrace{0...0}_{m_{i-1}}21\underbrace{9...9}_{l_i}78\underbrace{0...0}_{m_{i-1}}21\underbrace{9...9}_{l_{i-1}}78...\underbrace{0...0}_{m_1}21\underbrace{9...9}_{l_1}78$$

with $m_i$ and $l_i$ being 0 or positive integers then $M$ is *(10, 4)*- reverse multiple. Here,

$$i = \lceil\frac{number\ of\ occurence\ of\ 219\cdots978}{2}\rceil.$$

**Lemma 2.** If $M$ is an integer of the form

$$10\underbrace{9...9}_{l_1}89\underbrace{0...0}_{m_1}\cdots\underbrace{0...0}_{m_{i-1}}10\underbrace{9...9}_{l_i}89\underbrace{0...0}_{m_i}10\underbrace{9...9}_{l_i}89\underbrace{0...0}_{m_{i-1}}\cdots\underbrace{0...0}_{m_1}10\underbrace{9...9}_{l_1}89$$

or

$$10\underbrace{9...9}_{l_1}89\underbrace{0...0}_{m_1}\cdots 10\underbrace{9...9}_{l_{i-1}}89\underbrace{0...0}_{m_{i-1}}10\underbrace{9...9}_{l_i}89\underbrace{0...0}_{m_{i-1}}10\underbrace{9...9}_{l_{i-1}}89...\underbrace{0...0}_{m_1}10\underbrace{9...9}_{l_1}89$$

with $m_i$ and $l_i$ being 0 or positive integers then $M$ is *(10, 9)*- reverse multiple. Here,

$$i = \lceil\frac{number\ of\ occurence\ of\ 109\cdots989}{2}\rceil.$$

**Theorem 2.** An integer is *(10, k)*-reverse multiple if and only if k is 1, 4 or 9.

**Proof.** Using the previous Theorem and Lemmas, it suffices to prove that k is not equal to neither any of 2, 3, 5, 6, 7 nor 8.
We have 6 cases for k, in each case, we may set $M$ in digit form; i.e.
$M = a_n a_{n-1} a_{n-2} \cdots a_2 a_1 a_0$. It is clear that $a_n$ and $a_0$ are non-zero digits.
Define $c_i$ to be the left digit when multiplying k by $a_{n-i-1}$. These forms for $M$ and $c_i$ will be used throughout this paper. We prove all the cases by contradiction.

*Case k = 2.*
$M$ is *(10, 2)*- reverse multiple then $a_n \equiv 2a_0 (mod 10)$, $2a_n \leq a_0 \leq 9$ and 2 divides $a_n$.
This implies that $a_n = 2$ or 4.
If $a_n = 4$ then $2a_0 \equiv 4(mod 10)$ which contradicts the fact that $a_0$ is either equal to 8 or to 9.
If $a_n = 2$ then $a_0 \geq 4$ and $2a_0 \equiv 2(mod 10)$. The latter implies that $a_0 = 6$.
Our number now has the form $2a_{n-1}a_{n-2}a_{n-3}\cdots a_3 a_2 a_1 6$. It is easy to see that $2a_1 + 1 \equiv a_{n-1}(mod 10)$ implying that $a_{n-1} - 1$ is even, and $2a_{n-1} + c_1 = a_1 + 20$ where $c_1$ is the left digit when multiplying 2 by $a_{n-2}$. Having $c_1 - a_1$ a positive multiple of 2 that is less than 10, we get that $c_1 - a_1 = 2, 4, 6$ or 8.

If $c_1 - a_1 = 2$ then $a_{n-1} = 9$ and $a_1 = 4$ or $9$. It is clear that $a_1 \neq 9$ else $c_1 = 11$. Thus, we get that $a_1 = 4$ and $c_1 = 6$ resulting in $2a_{n-2} + c_2 = a_2 + 60$. The latter has no digit solution.

If $c_1 - a_1 = 4$ or $c_1 - a_1 = 8$ then $a_{n-1} = 8$ and $a_{n-1} = 6$ respectively which is contradiction to the fact that $a_{n-1}$ is odd.

If $c_1 - a_1 = 6$ then $a_{n-1} = 7$ and $a_1 = 3$ or $8$ and both give contradiction. The case $a_1 = 3$ gives $c_1 = 9$ and thus equation $2a_{n-2} + c_2 = 30 + a_2$ has no digit solution. And the case $a_1 = 8$ implies $c_1 = 11$.

*Case k = 3.*
M is *(10, 3)*- reverse multiple then $a_n \equiv 3a_0 \pmod{10}$ and $3a_n \leq a_0 \leq 9$. This implies that $a_n = 1$, 2 or 3. The value of $a_n$ is neither 2 nor 3 else $a_0$ is 4 and 1 respectively which contradicts the fact that $a_0 \geq 6$. Moreover, $a_n \neq 1$ since $a_0 = 7$ and the equation $3a_{n-1} + c_1 = a_1 + 40$ has no digit solution.

*Case k = 5.*
M is *(10, 5)*- reverse multiple then $a_n \equiv 5a_0 \pmod{10}$, $5a_n \leq a_0 \leq 9$ and $a_n$ is a multiple
of 5. This implies that $a_n = 5$ which is impossible since $5a_n \leq a_0$ which is a digit.

*Case k = 6.*
M is *(10, 6)*- reverse multiple then $a_n \equiv 6a_0 \pmod{10}$, $6a_n \leq a_0 \leq 9$ and 2 divides $a_n$. It is clear that $a_n \geq 2$ which is impossible since $a_0 \geq 6a_n \geq 12$ is a digit.

*Case k = 7.*
M is *(10, 7)*- reverse multiple then $a_n \equiv 7a_0 \pmod{10}$ and $7a_n \leq a_0 \leq 9$. This implies that
$a_n = 1$ and $a_0 = 3$ which contradicts the fact that $7a_n \leq a_0$.

*Case k = 8.*
M is *(10, 8)*- reverse multiple then $a_n \equiv 8a_0 \pmod{10}$, $8a_n \leq a_0 \leq 9$ and 2 divides $a_n$.
This implies that $a_n = 1$ which is impossible since $8a_0 \equiv 1 \pmod{10}$ is not solvable in $Z_{10}$.

## 3. Necessary conditions

In this section, we find necessary condition for *(10, k)*- reverse multiples by presenting our next theorems.

**Theorem 3.** If $M$ is *(10, 4)*- reverse multiple then $M$ has the form given in Lemma 1.

**Proof.** M is *(10, 4)*- reverse multiple then $a_n \equiv 4a_0 \pmod{10}$, $4a_n \leq a_0 \leq 9$ and 2 divides $a_n$. By solving these equations, we get that $a_n = 2$ and $a_0 = 8$.
Our number now is $2a_{n-1}a_{n-2}a_{n-3} \cdots a_3 a_2 a_1 8$ with $a_{n-1}$ being odd as 4 is a divisor of $(a_{n-1}2)$. It is easy to see that $4a_1 + 3 \equiv a_{n-1} \pmod{10}$ and that $4a_{n-1} + c_1 = a_1$. The digit $a_1 - c_1$ is a non zero positive multiple of 4 which implies that $a_1 - c_1 = 4$ or 8. But $a_1 - c_1 \neq 8$ else $a_{n-1} = 2$ which is not odd. Thus, $a_1 - c_1 = 4$.
Substituting $a_1 - c_1 = 4$ in our equations, we get that $a_{n-1} = 1$ and $4a_1 \equiv 1 \ (mod\ 10)$. The latter implies that $a_1 = 2$ or $7$, $a_1 \neq 2$ else $c_1 = -2$.
Thus $a_1 = 7$ and $c_1 = 3$.
Our number now is $21a_{n-2}a_{n-3}\cdots a_3a_278$ and our new equations are
$$4a_2 + 3 \equiv a_{n-2}(mod\ 10)$$
$$4a_{n-2} + c_2 = 30 + a_2.$$
We get that $a_{n-2}$ is odd and 4 is a divisor of $a_2 - c_2 + 2$; i.e $a_2 - c_2 + 2 = 0, \pm 4$ or $\pm 8$. The cases $a_2 - c_2 + 2 = 4, -4, -8$ imply $a_{n-2} = 8, 6$ and $5$. The first two of them contradict the fact that $a_{n-2}$ is odd and the third contradicts the fact that $c_2 = a_2 + 10$ is digit.
We are left with the following cases

*Case* $a_2 - c_2 + 2 = 8$.

In this case, we get that $a_{n-2} = 9$ and $a_2 = 4$ or $9$. But $a_2 \neq 4$ else $6 - c_2 = 8$.
Thus, $a_2 = 9$ and $c_2 = 3$.
Our number now is $219a_{n-3}a_{n-4}\cdots a_4a_3978$ and our new equations are
$$4a_3 + 3 \equiv a_{n-3}(mod\ 10)$$
$$4a_{n-3} + c_3 = 30 + a_3.$$
The above equations have the same solution as that for $a_2$ and $a_{n-2}$. So, by repeating the same procedure we get that $a_3 - c_3 + 2 = 8$ or $a_3 - c_3 + 2 = 0$.
Continuing on this pattern, we get $a_k - c_k + 2 = 8$ or $a_k - c_k + 2 = 0$. And if $a_k - c_k + 2 = 8$, we get $a_k = a_{n-k} = 9$ and $c_k = 3$. Thus, $M = 219\cdots 978$.

The case $a_k - c_k + 2 = 0$ will be resolved in Case $a_2 - c_2 + 2 = 0$.

*Case* $a_2 - c_2 + 2 = 0$.
We have $c_2 = a_2 + 2 \leq 9$ which implies that $a_2 \leq 7$ and $a_{n-2} = 7$. And $4a_2 \equiv 4(mod\ 10)$ implies that $a_2 = 1$ or $6$. But $a_2 \neq 6$ else $c_2 = 8$ and the equation $4a_{n-3} + c_3 = 80 + a_3$ has no digit solution. Thus, $a_2 = 1$, $a_{n-2} = 7$ and $c_2 = 3$.
Our number now is $217a_{n-3}a_{n-4}\cdots a_4a_3178$ and our new equations are
$$4a_3 \equiv a_{n-3}(mod\ 10)$$
$$4a_{n-3} + c_3 = 30 + a_3.$$
It is easy to see that the solution of our new equations is $a_3 = 2$, $a_{n-3} = 8$ and $c_3 = 0$.
Our number now is $2178a_{n-4}a_{n-5}\cdots a_5a_42178$ and our new equations are
$$4a_4 \equiv a_{n-4}(mod\ 10)$$
$$4a_{n-4} + c_4 = a_4 \leq 9.$$
The above equations imply that $a_{n-4}$ is even and less than or equal to 2 and that $a_4 - c_4 = 0, 4$ or $8$. Thus, $a_{n-4} = 0$ or $2$.
If $a_{n-4} = 2$ then $a_4 = 3$ or $8$ and $a_{n-4} \neq 3$ since $8 + c_4 = 3$ has no digit solution.
Having $a_4 = 8$, we get that $c_4 = 0$ and by repeating the same procedure as done before, we get that $a_{n-5} = 1$ then $a_5 = 7$.
If $a_{n-4} = 0$ then $a_4 = 0$ or $5$. But $a_4 \neq 5$ else $c_4 = 5$ and $4a_{n-5} + c_5 = 50 + a_5$ which has no digit solution. Thus, $a_4 = c_4 = 0$.
Our number now is $21780a_{n-5}a_{n-6}\cdots a_6a_502178$ and our new equations are
$$4a_5 \equiv a_{n-5}(mod\ 10)$$
$$4a_{n-5} + c_5 = a_5 \leq 9.$$
We get the same equation as that for $a_0$ and $a_n$ with a slight difference that $a_5$ and $a_{n-5}$ can be equal to zero. Using the results of $a_0$ and $a_n$, we get that $a_5 = 8$ and $a_{n-5} = 2$ or $a_5 = a_{n-5} = 0$.
Continuing on this pattern, we get the form of $M$ that is stated in Theorem 3.

Next, we present a theorem regarding the *(10, 9)*- reverse multiples. The proof of it is the same as that of the previous theorem but with less and easier computations.

**Theorem 4.** If $M$ is *(10, 9)*- reverse multiple then $M$ has the form of Lemma 2. Finally, we present our main theorem which is a result of our proved lemmas and theorems in the paper.

**Theorem 5.** An integer $M$ is *(10, k)*- reverse multiple if and only if one of the following conditions is satisfied

1. $k = 1$ and $M$ is a palindrome;

2. $k = 4$ and $M$ has the forms

$$21\underbrace{9\ldots9}_{l_1}78\underbrace{0\ldots0}_{m_1}\cdots\underbrace{0\ldots0}_{m_{i-1}}21\underbrace{9\ldots9}_{l_i}78\underbrace{0\ldots0}_{m_i}21\underbrace{9\ldots9}_{l_i}78\underbrace{0\ldots0}_{m_{i-1}}\ldots\underbrace{0\ldots0}_{m_1}21\underbrace{9\ldots9}_{l_1}78$$

or

$$21\underbrace{9\ldots9}_{l_1}78\underbrace{0\ldots0}_{m_1}\cdots 21\underbrace{9\ldots9}_{l_{i-1}}78\underbrace{0\ldots0}_{m_{i-1}}21\underbrace{9\ldots9}_{l_i}78\underbrace{0\ldots0}_{m_{i-1}}21\underbrace{9\ldots9}_{l_{i-1}}78\ldots\underbrace{0\ldots0}_{m_1}21\underbrace{9\ldots9}_{l_1}78$$

with $m_i$ and $l_i$ being 0 or positive integers, $i = \lceil\frac{number\ of\ occurence\ of\ 219\cdots978}{2}\rceil$.

3. $k = 9$ and $M$ has the forms

$$10\underbrace{9\ldots9}_{l_1}89\underbrace{0\ldots0}_{m_1}\cdots\underbrace{0\ldots0}_{m_{i-1}}10\underbrace{9\ldots9}_{l_i}89\underbrace{0\ldots0}_{m_i}10\underbrace{9\ldots9}_{l_i}89\underbrace{0\ldots0}_{m_{i-1}}\ldots\underbrace{0\ldots0}_{m_1}10\underbrace{9\ldots9}_{l_1}89$$

or

$$10\underbrace{9\ldots9}_{l_1}89\underbrace{0\ldots0}_{m_1}\cdots 10\underbrace{9\ldots9}_{l_{i-1}}89\underbrace{0\ldots0}_{m_{i-1}}10\underbrace{9\ldots9}_{l_i}89\underbrace{0\ldots0}_{m_{i-1}}10\underbrace{9\ldots9}_{l_{i-1}}89\ldots\underbrace{0\ldots0}_{m_1}10\underbrace{9\ldots9}_{l_1}89$$

with $m_i$ and $l_i$ being 0 or positive integers, $i = \lceil\frac{number\ of\ occurence\ of\ 109\cdots989}{2}\rceil$.

**Notes:**
1. There are infinitely many *(10, 4)* and *(10, 9)*- reverse multiples.
2. For a fixed number of digits, there are finite number of *(10, 4)* and *(10, 9)*- reverse multiples. For example, there are two (10,4)- reverse multiples of length equals 8.
3. The first *(10, 4)* - reverse multiples are 2178, 21978, 219978, 2199978, 21782178, 21999978, 217802178 .

The following table represents an arbitrary list of reverse multiples.

| number | $l_1$ | $m_1$ | $l_2$ | $m_2$ | type |
|---|---|---|---|---|---|
| 21782178 | 0 | 0 | | | (10, 4) |
| 21782197800219782178 | 0 | 0 | 1 | 2 | (10, 4) |
| 1089010999890 1089 | 0 | 1 | 3 | | (10, 9) |

DEPARTMENT OF MATHEMATICS AND PHYSICS, LEBANESE INTERNATIONAL UNIVERSITY, LEBANON
Email address: madeline.tahan@liu.edu.lb